\documentclass[12pt]{article}
\usepackage{amsmath,amssymb,amsthm,amsfonts}
\usepackage{latexsym}
\begin{document}

\newtheorem{prop}{Proposition}[section]
\newtheorem{cor}{Corollary}[section] 
\newtheorem{theo}{Theorem}[section]
\newtheorem{lem}{Lemma}[section]
\newtheorem{rem}{Remark}[section]
\newtheorem{con}{Conjecture}[section]
\newtheorem{as}{Assumption}[section]
\newtheorem{de}{Definition}[section]

\setcounter{page}{1}
\renewcommand{\theequation}{\thesection.\arabic{equation}}
\begin{center}
{\Large \bf Asymptotics for the size of the largest
component scaled to "log n" in inhomogeneous random graphs.}
\end{center}

\begin{center}
TATYANA S. TUROVA\footnote{Research was 
supported by the Swedish Natural Science Research
Council.

2000 {\it  Mathematics Subject Classification}: 60C05; 05C80.

{\it Key words and phrases:} inhomogeneous random graphs, largest
connected component, multi-type branching processes, phase transition.
} 
\end{center}

\begin{center}
{\it Mathematical Center, University of
Lund, Box 118, Lund S-221 00, 
Sweden. }
\end{center}

\begin{abstract}
We study the inhomogeneous random graphs in the subcritical case.
We derive an  exact formula for the 
size of the largest connected component scaled to $\log n$ where $n$ is the size of the graph.
This generalizes the recent result for the "rank 1 case". Here we discover that 
the same well-known equation for the survival probability, 
whose positive solution determines the 
asymptotics of the 
size of the largest component in the supercritical case, 
plays the crucial role in the subcritical case as well. 
But now these are the {\it negative} solutions which come into play.
\end{abstract}

\medskip

\renewcommand{\theequation}{\thesection.\arabic{equation}}
\section{Introduction.}
\setcounter{equation}{0}

\subsection{Inhomogeneous random graphs.}

A general inhomogeneous
random graph model which comprises numerous previously known models,
was introduced and studied in great details by Bollob{\'a}s,   Janson and   Riordan 
\cite{BJR}. Let us recall the basic definition of the inhomogeneous
random graph
$G^{\cal V}(n,\kappa )$ 
with a vertex space
$${\cal V}=(S,\mu, (x_1^{(n)}, \ldots, x_{n}^{(n)})_{n\geq 1}).$$
Here $S$ is a separable metric space and $\mu$ is a Borel probability
measure on $S$.
No
relationship is assumed between $x_i^{(n)}$ and $x_i^{(n')}$, but to
simplify notations we
shall write further $(x_1, \ldots, x_{n})=(x_1^{(n)}, \ldots, x_{n}^{(n)})$.
For each $n$ let
$(x_1, \ldots, x_{n})$ be
a deterministic or random
sequence  of points in $S$, such that for any
$\mu$-continuity set $A\subseteq S$
\begin{equation}\label{set}
    \frac{\#\{i: x_i\in A\}}{n}\stackrel{P}{\rightarrow}   \mu (A).
\end{equation}
Given the sequence $x_1, \ldots, x_{n}$, we let  $G^{\cal V}(n,\kappa)$
 be the random graph on $\{1, \ldots, n\}$, such that any
two vertices $i$ and $j$ are connected by an edge independently of the others
and with a probability
\begin{equation}\label{pe}
p_{x_i x_j}(n)= \min\{\kappa (x_i,x_j)/n,1\},
\end{equation}
where  the kernel
$\kappa$ is a symmetric non-negative measurable function on $S\times
S$.
We shall also assume that kernel
$\kappa$ is {\it graphical} on $\cal V$, which means that

(i) $\kappa$ is continuous $a.e.$ on $S\times S$;

(ii) $ \kappa \in L^1(S \times S,  \mu\times \mu)$;

(iii)
\[\frac{1}{n} {\bf E} e\Big(G^{\cal V}(n,\kappa )\Big)
\rightarrow
\frac{1}{2}\int_{S^2}\kappa (x,y)d {\mu}(x)d
{\mu}(y),\]
where $e(G)$ denotes the number of edges in a graph $G$.

Denote  $ C_1 \Big( G \Big)$  the size of
the largest
connected component in a graph $G$. This
is the most studied characteristic of the random graphs.
In particular, the famous phenomena of phase transition is seen in 
the abrupt change of the value  $ C_1 \Big( G \Big)$
 depending on the parameters of the model.
It appears that there is a close connection between 
$ C_1 \Big( G \Big)$ and a survival probability of  a certain  multi-type 
 Galton-Watson process $B_{\kappa}(x)$ defined below.

\noindent
\begin{de}\label{d1}
 The type space of  $B_{\kappa}(x)$ is $S$, and initially there
 is a single particle of type $x \in S$. Then 
 at any step,
 a particle of type $x \in S$ is replaced in the next
 generation by a set of particles where the number of  particles of
 type $y$ has a Poisson distribution
with intensity 
 $\kappa
 (x,y)d {\mu}(y).$
\end{de}
\noindent
Let $\rho_{\kappa}(x)$ denote the survival probability of $B_{\kappa}(x)$. 
Then Theorem 3.1 from \cite{BJR} states   that
\begin{equation}\label{B}
\frac{C_1 \Big(G^{\cal V}(n,\kappa)
    \Big)}{ n }\stackrel{P}{\rightarrow} \rho_{\kappa}: = \int_{S} \rho_{\kappa}(x)
  d\mu(x).
  \end{equation}
It was also proved in \cite{BJR} that $\rho_{\kappa}(x)$
 is the maximum
  solution to
\begin{equation}\label{+}
f(x) = 1- e^{-T_{\kappa}[f](x)},
\end{equation}
where
the
integral operator $T_{\kappa}$  is defined by
\begin{equation}\label{Tk}
(T_{\kappa}f)(x)=\int_S{\kappa}(x,y)f(y)d\mu(y).
\end{equation}
Whether $\rho_{\kappa}$ is zero or strictly positive depends only on
the norm 
\[\| T_{\kappa}\|= \sup\{\| T_{\kappa}f\|_2:f\geq 0, \|f\|_2\leq1\}.\]
Namely, due to Theorem 3.1 from \cite{BJR} 
\begin{equation}\label{rok}
\rho_{\kappa} \ 
\left\{ \
\begin{array}{ll}
> 0, & \mbox{ if } \ \| T_{\kappa}\|>1, \\
= 0, & \mbox{ if } \ \| T_{\kappa}\|\leq 1.
\end{array}
\right.
\end{equation}
This  together with (\ref{B}) tells us that in the subcritical case, i.e., when $\| T_{\kappa}\| \leq  1$, we have
$C_1\Big(G^{\cal V}(n,\kappa )\Big)=o_P(n).$ 

Under an additional  assumption 
$$
\sup_{x,y}\kappa (x,y)<\infty $$
Theorem 3.12 in \cite{BJR} proves
 that if $\| T_{\kappa}\| < 1$ then
 $C_1\left( G^{\cal V}(n,\kappa )
\right)=O(\log n)$  with probability tending to one as $n \rightarrow
 \infty$. 

On the other hand, as it was pointed out in \cite{BJR},   when 
the kernel is unbounded, the condition $\| T_{\kappa}\| < 1$
is not sufficient for  the size of the largest component  
to be of order $\log n$; an example is the model of random growth from \cite{BJR1}. 
Recently Janson showed in  \cite{J1} that a  subcritical inhomogeneous random graph
can also have the largest component of order $n^{1/\gamma}$
under the assumption of a power law degree distribution with exponent $\gamma+1$, $\gamma>2$.

Here we describe  sufficient conditions under which  $C_1\left( G^{\cal V}(n,\kappa )
\right)/\log n $ converges in probability to a finite  constant
even for unbounded kernels.
The exact value of this constant
untill recently
was known only for the
Erd{\" o}s-R{\' e}nyi random graph
\cite{ER}. The first related result for the inhomogeneous model 
but only in the rank 1 case, i.e., when 
\begin{equation}\label{k1}
{\kappa}(x,y)=\Phi(x)\Phi(y),
\end{equation}
was derived  in \cite{Tlog}. 
However, in \cite{Tlog} the formula for the
asymptotics  of $C_1\left( G^{\cal V}(n,\kappa )
\right)/ \log n$ is given in terms of 
function $\Phi$ and thus  is not applicable for a general kernel. Here we shall
consider a more general situation, which includes as well case 
(\ref{k1}).

\subsection{Results.}
Denote ${\cal X}(x)$ the size of the  total progeny of $B_{\kappa}(x)$, 
and let
\begin{equation}\label{1}
r_{\kappa}=sup \, \{z\geq 1: \int_S  {\bf E}\left(
  z^{ {\cal X}(x)}\right)
d\mu(x)\  <\infty\}.
\end{equation}
We will show that $r_{\kappa}$ is the determining
parameter for the
size $C_1 \Big(
G^{\cal V}(n,\kappa )
    \Big)$ in the subcritical case. 
In particular, we need to know whether
$r_{\kappa}=1$ or  $r_{\kappa}>1$.
    Therefore first we shall study
$r_{\kappa}$. One should notice the direct relation of $r_{\kappa}$ to
the tail of distribution of the total progeny  ${\cal X}(x)$.
In particular, if the tail of distribution of ${\cal X}(x)$ decays
exponentially, $r_{\kappa}$ defines the constant in the exponent. In
the case of a single-type branching process the exact result on the
relation between $r_{\kappa}$ and the distribution of the total
progeny was proved in \cite{O}.

Note that when 
\begin{equation}\label{R*}
\int_S \kappa(x,y)d\mu(y)<\infty \ \mbox{ for all} \ x\in S,
\end{equation}
 Lemma 5.16 in \cite{BJR} states the following: if $\| T_{\kappa}\|>1$
then
  $\rho_{\kappa}>0 $ on a set of positive measure. 
  This means that ${\cal X}=\infty$ on a set of positive measure, which  immediately implies 
\begin{equation}\label{rRev}
r_{\kappa}=1, \ \mbox{ if } \ \| T_{\kappa}\|>1.
\end{equation}
 
We shall assume from now on that
\begin{equation}\label{inf}
\inf _{x,y \in S}{\kappa}(x,y)>0.
\end{equation}

\begin{theo}\label{T1}
$r_{\kappa}$ is the supremum value of all
$z\geq 1$ for which equation
\begin{equation}\label{-}
g(x)=ze^{T_{\kappa}[g-1](x)}
\end{equation}
has  $a.s.$ (i.e., $\mu-a.s.$) finite solution  $g\geq 1$.
\end{theo}

Theorem \ref{T1}
 yields immediately the following criteria.
\begin{cor}\label{C3}
 $r_{\kappa} >1$ if and only if at least for some $z>1$ equation 
(\ref{-})
has an $a.s.$ finite solution  $g>1$. Otherwise,  $r_{\kappa} =1$. \hfill$\Box$
\end{cor}

Next we  extend statement (\ref{rRev}) for the case $\| T_{\kappa}\|=1$.
\begin{cor}\label{C4} Let $\kappa$ satisfy (\ref{R*}). Then
\begin{equation}\label{r=1}
r_{\kappa} =1 \ \  \mbox{ if \ }  \ \  \| T_{\kappa}\| \geq 1.
\end{equation}
\end{cor}

Theorem \ref{T1} and Corollary \ref{C3} will allow us to derive some sufficient conditions for $r_{\kappa}>1$.
Let $T_{\kappa}$ have a finite Hilbert-Schmidt norm,
i.e., 
\begin{equation}\label{HS}
\| T_{\kappa}\|_{HS}:= \|\kappa \|_{L^2(S\times S)}=
\left(\int_S\int_S{\kappa}^2(x,y)d\mu(x)d\mu(y)\right)^{1/2}<\infty.
\end{equation}
Define  
\begin{equation}\label{dpsi}
\psi(x)=\left(\int_S{\kappa}^2(x,y)d\mu(y)\right)^{1/2},
\end{equation}
and assume  that for some positive constant $a>0$ 
\begin{equation}\label{As1}
\int_{S}e^{a \psi(x)} d\mu(x)< \infty.
\end{equation}
\begin{theo}\label{C1} 
Let $\kappa$ satisfy (\ref{As1}).
Then
\begin{equation}\label{r>1}
r_{\kappa} >1 \ \  \mbox{ if \ }  \ \  \| T_{\kappa}\| < 1
\end{equation}
and at least one of the following conditions is satisfied

$(C1)$ $ \sup_{x,y \in S} {\kappa}(x,y)<\infty$, or

$(C2)$ $\| T_{\kappa}\|_{HS}<1$, or

$(C3)$ $S\subseteq {\bf R}$ and $\kappa(x,y) $ is 
   non-decreasing in both arguments, and such that 
   for some constant $c_1>0$
\begin{equation}\label{TT}
   \kappa(x,y) \leq c_1 T_{\kappa}[1](x) T_{\kappa}[1](y),
\end{equation}
for all $x,y\in S$.

\end{theo}

\begin{rem} Condition $S\subseteq {\bf R}$ in $(C3)$
one can  replace by a condition that  space $S$ can be partially
ordered.
 \end{rem}

Observe that for all kernels  
\[\| T_{\kappa}\| \leq \| T_{\kappa}\|_{HS} ,\]
where equality holds only in the rank 1 case (\ref{k1}). Hence, under assumption (\ref{As1})
in the
rank 1  case condition  $\| T_{\kappa}\|<1$ is
sufficient and necessary for $r_{\kappa}>1$.
\bigskip

Consider now model $G^{\cal V}(n,\kappa )$ which satisfies 
(\ref{set}) and (\ref{pe}).
To be able to approximate a component in $G^{\cal V}(n,\kappa )$ by
a branching process  we need some additional conditions on the distribution
of the types of vertices $x_1, \ldots, x_n$.

\begin{as}\label{Ass1} Let
$S \subseteq \{1, 2, \ldots \}$ be finite or countable, and 
let
for any $\varepsilon>0$ and $q>0$
\begin{equation}\label{As3}
{\bf P} \left\{
   \frac{\#\{1\leq i\leq n: x_i=x\}}{n} - \mu(x)  \leq
  \varepsilon e^{qT_{\kappa}[1](x)}\mu(x) , \mbox{ for all }  x\in S \right\}
 \rightarrow 1 
\end{equation}
as $ n\rightarrow \infty$.
\end{as}

\noindent
Notice, that in the case when $S$ is finite convergence (\ref{As3})
follows simply by the assumption (\ref{set}). Some examples of
the models with countable $S$ which satisfy condition  (\ref{As3})
one can find in
\cite{TV}.

\noindent
\begin{theo}\label{T1A} 
Let $\kappa$ satisfy
(\ref{As1}), as well as at least one of the 
  conditions 
$(C1)$ or
$(C3)$ from Theorem \ref{C1}.
Then under Assumption \ref{Ass1}
\begin{equation}\label{lp}
\frac{C_1 \Big(G^{\cal V}(n,\kappa )
    \Big)}{\log n }
\stackrel{P}{\rightarrow} \frac{1}{\log r_{\kappa}},
\end{equation}
where
\begin{equation}\label{r1}
r_{\kappa} \left\{  \begin{array}{ll}
 >1, & \mbox{ if \ } \| T_{\kappa}\|<1  , \\ \\
=1, & \mbox{ if \ }   \| T_{\kappa}\| \geq 1.
\end{array}
\right.
\end{equation}
\end{theo}

\bigskip

 Theorem \ref{T1A}
 provides
   sufficient conditions when 
convergence
(\ref{lp}) 
 takes place even for unbounded kernels. 
Observe, however,  that condition (\ref{As1}) seems to be  necessary as well.
In particular, in the "rank 1" case (\ref{k1}) condition (\ref{As1})
excludes possibility of a power law degree distribution. The later is
proved (\cite{J1}) to  yield order $n^{1/\gamma}$ ($\gamma>2$) for
the largest component in a subcritical graph. 

Clearly, Theorem \ref{T1A} complements statement (\ref{B}) together
with (\ref{rok}). There is even a direct relation between the values
$r_{\kappa}$ and $\rho_{\kappa}$ as we shall see now.
Setting $f(x)=-(g(x)-1)$ in (\ref{-}),  we get from Corollary \ref{C3}
that
 $r_{\kappa} >1$ {\it if and only if at least for some $z>1$ equation 
\begin{equation}\label{1-}
f(x)=1-ze^{-T_{\kappa}[f](x)}
\end{equation}
 has an $a.s.$ finite solution  $f<0$.}
Notice that when $z=1$ equation (\ref{1-}) coincides with (\ref{+}).
This observation leads to  a surprising direct relation to the
supercritical case.

\begin{prop}\label{Prev}
Let $\kappa$ satisfy  (\ref{inf}). Then 
 $r_{\kappa} >1$ if  equation (\ref{+})  has an $a.s.$ finite solution
$f$ such that  $\sup_{x\in S} f(x)<0$.
\end{prop}

In the case of 
a homogeneous Erd{\" o}s-R{\' e}nyi graph
$G_{n,p}$ (consult, e.g., \cite{B})
where the probability of any edge is $p=c/n$,
the relation between the supercritical and subcritical cases is most
transparent. Placing $G_{n,p}$ into the general definition of an 
inhomogeneous random graph model
gives us
 $|S|=1$ and  $\kappa \equiv c$.
The corresponding branching process $B_{\kappa}$ (see Definition \ref{d1}) has $Po(c)$
distribution of the offspring (of a single type). The  survival
probability  $\rho_c$ of this process is again the maximum solution
to (\ref{+}), which takes a simple form
\begin{equation}\label{1+a}
f=1-e^{-cf}.
\end{equation}
By Corollary \ref{C3} we have here $r_c>1$ if and only if equation (\ref{-}), which in
this case is
\begin{equation}\label{-a}
g=ze^{c(g-1)},
\end{equation}
has a finite solution $g>1$ for some $z>1$.
It is straightforward to check that (\ref{-a}) has a finite positive solution 
for some $z>1$ if and only if equation (\ref{1+a}) has a strictly negative
solution (or, else, if and only if $c<1$).

Therefore we may conjecture
that the condition in 
Proposition \ref{Prev} is necessary as well.

\begin{con}\label{Crev}
Let $\kappa$ satisfy  (\ref{inf}). Then
 $r_{\kappa} >1$ if  and only if equation (\ref{+})  has an $a.s.$ finite solution
$f<0$.
\end{con}

Observe,  that while all the nonnegative solutions to
(\ref{+}) are bounded by 1, the nonpositive ones can be
unbounded. This certainly makes a  difference for the analysis. To
surpass the difficulties we introduced condition $(C3)$, which 
resembles  a rank 1 case. One may believe that the results of Theorem
\ref{C1}
and Theorem
\ref{T1A} should hold in a much more general situation than we are able to
treat  here.

In the special
"rank 1 case" (\ref{k1}) convergence (\ref{lp}) was previously established in
\cite{Tlog}
under some additional conditions on function
$\Phi$. 
Note that  in the rank 1 case (see for the details \cite{Tlog}) one can derive an explicite  formula for
$r_{\kappa}$.
Clearly, condition (\ref{k1})  implies (\ref{TT}). On the other hand,
there are kernels which satisfy the condition
(\ref{TT}) but not (\ref{k1}). These are for example,
$ {\kappa}(x,y) =\widetilde{\kappa}(x\vee y),$
where $\widetilde{\kappa}$ is a positive monotone increasing
function on $S$, such that $\int_S\widetilde{\kappa}^2d\mu < \infty$. Models with kernels of this type were considered, e.g.,
in \cite{BJR} and \cite{T}.

\bigskip
\section{The generating function for the progeny of a branching
  process.}
\setcounter{equation}{0}

Recall that we denote
${\cal X}(x)$ the size of the  total progeny of $B_{\kappa}(x)$,
(see Definition \ref{d1} in the Introduction). 
We shall study  function
\[h_z(x)={\bf E} z^{{\cal X}(x)}, \ \ x\in S,\]
for $z\geq 1$ . 
It is standard to derive (consult, e.g., \cite{M}, Chapter 6)  that $h_z(k)$ as a generating function for a  branching
process
satisfies the following equation
\[h_z(k) =z \exp {\left\{
\int_{S} \kappa  (k,x)(h_z(x)-1 ) d \mu(x) 
\right\}}, \ \ \ \ k\in S,\]
or in a functional form:
\begin{equation}\label{A1}
h_z = ze^{ T_{\kappa}[h_z-{1}]}=:\Phi_{z,\kappa}h_z.
\end{equation}

\begin{theo}\label{T3}
For any $z\geq 1$ function $h_z$ is the minimal solution $f\geq 1$
 to the equation
 \begin{equation}\label{A3}
f=\Phi_{z,\kappa}f,
\end{equation}
and moreover
\begin{equation}\label{Jn}
h_{z}=\lim_{k\rightarrow
  \infty}\Phi_{z,\kappa}^k[1].
\end{equation}
\end{theo}
\noindent
{\bf Proof.}
Let us denote ${\cal X}_k(x)$, $k\geq 0$, the number of the offspring
of the process $B_{\kappa}(x)$ in the first $k$ generations.
In particular, ${\cal X}_0(x)=1$ and
\[{\cal X}_1(x)=_d 1 + P_{x},\]
where $P_{x}$ is the number of the offspring of a particle of type
$x$, among which the number
of particles of each
type $y\in S$ has $Po(\kappa(x,y)d\mu(y))$-distribution.
Let
\[h_{k,z}(x)= {\bf E}z^{{\cal
    X}_k(x)}\]
for $k\geq 0$. It is straightforward  to derive that 
\[h_{1,z}(x)= {\bf E}z^{{\cal
    X}_1(x)}=\Phi_{z,\kappa}[z](x)=\Phi_{z,\kappa}[h_{0,z}](x), \]
and similarly for any $k\geq 1$
\[h_{k+1,z}(x)=\Phi_{z,\kappa}[h_{k,z}](x).\]
Noticing that $h_{0,z}(x)=z=\Phi_{z,\kappa}[1](x)$ for all $x\in S$, we derive from here
\begin{equation}\label{A9}
h_{k,z}(x)=\Phi_{z,\kappa}^{k+1}[1](x).
\end{equation}
Obviously, $h_{k,z}(x) \nearrow h_{z}(x)$, i.e.,
\[h_{z}(x)=\lim_{k\rightarrow \infty}\Phi_{z,\kappa}^k[1](x) \]
for all $x\in S$.
By the monotone convergence 
\[T_{\kappa}[h_{z}](x)= \int_{S} \kappa  (x,y)\lim_{k\rightarrow
  \infty}\Phi_{z,\kappa}^k[1](y) d \mu(y) 
=
\lim_{k\rightarrow \infty}T_{\kappa}[\Phi_{z,\kappa}^k[1]](x),
\]
and therefore
\begin{equation}\label{A10}
\Phi_{z,\kappa}[h_{z}](x)=e^{\lim_{k\rightarrow \infty}
    T_{\kappa}[\Phi_{z,\kappa}^k[1]-1](x)}= \lim_{k\rightarrow \infty}
\Phi_{z,\kappa}[\Phi_{z,\kappa}^k[1]](x) =h_{z}(x).
 \end{equation}
Hence,  $h_{z}=\lim_{k\rightarrow
  \infty}\Phi_{z,\kappa}^k[1]$ is a solution to 
(\ref{A3}).

Since $\Phi_{z,\kappa}$ is monotone and
$\Phi_{z,\kappa}[1](x)=z\geq 1,$
it follows by induction 
that
\begin{equation}\label{O10}
h_{z}(x)=\lim_{k\rightarrow
  \infty}\Phi_{z,\kappa}^k[1](x)\geq 1
\end{equation}
for all $x\in S$.

Finally, we show that $h_{z}$ is the minimal solution $f\geq 1$ to 
(\ref{A3}).
Assume, that there is a solution $f\geq 1$ such that
$1\leq f(x)<h_{z}(x)$ for some $x$.
Then due to the monotonicity of $\Phi_{z,\kappa}$
 we have also
\[\Phi_{z,\kappa}^k[1](x) \leq \Phi_{z,\kappa}^k[f](x) = f(x)
< h_{z}(x)=\lim_{N\rightarrow
  \infty}\Phi_{z,\kappa}^N[1](x) \]
for all $k\geq 1$.
 Letting $k\rightarrow
  \infty$ in the last formula we come to the contradiction with
the strict inequality in the middle. Therefore $h_{z}$ is the minimal solution $f\geq 1$ to 
(\ref{A3}).
 \hfill$\Box$

\begin{rem}\label{R1}
If $f\geq 1$ satisfies (\ref{A3})
and $f(x)<\infty$ at least for some $x$, then
it follows straight from the definition of $\Phi_{z,\kappa}$
that $\int_{S} \kappa  (x,y)f(y)d \mu(y)<\infty$.
Hence, under assumption (\ref{inf}) if $f\geq 1$ satisfies (\ref{A3}) then either $f=\infty$ $a.s.$, or
$f<\infty$ $a.s.$, in which case also 
$\int_{S} \kappa  (x,y)f(y)d \mu(y) <\infty$ $a.s.$ The latter
together with the assumption (\ref{inf}) yields $f\in L_1(S, \mu)$ as well.
\end{rem}

\begin{rem}\label{RevC2} Theorem \ref{T1} and 
Corollary \ref{C3}
follow directly from  Theorem \ref{T3} and Remark \ref{R1}.
\end{rem}

Next we describe a sufficient condition when the minimal solution $f\geq 1$
 to the equation (\ref{A3}) is finite. 
\begin{lem}\label{L2}
  If $\Phi_{z,\kappa}f \leq f$ for some  $f\geq 1$,
  then there exists function $1\leq g\leq f$ which is a  solution to
  (\ref{A3}), i.e., 
  $\Phi_{z,\kappa}g=g.$
\end{lem}
\noindent
{\bf Proof.} (The proof almost repeats the one of Lemma 5.12 from
\cite{BJR}.)  The monotonicity of $\Phi_{z,\kappa}$  and assumption
$\Phi_{z,\kappa}f \leq f$ yield by induction
\[f\geq \Phi_{z,\kappa}f \geq \Phi_{z,\kappa}^2f\geq \ldots \ .\]
Since $f\geq 1$ we have for all $x$
\[\Phi_{z,\kappa}[f](x) =ze^{T_{\kappa}[f-1](x)} \geq z\geq 1,\]
which  implies by induction that also 
$ \Phi_{z,\kappa}^kf \geq 1 $
for all $k\geq 1$.
Hence the limit 
\[f(x)\geq g(x)=\lim_{k\rightarrow \infty}\Phi_{z,\kappa}^k[f](x) \geq 1\]
exists for every $x$. By the monotone convergence (repeat the
argument from (\ref{A10}))
$g$ is a solution to (\ref{A3}).
\hfill$\Box$

  \bigskip

\begin{theo}\label{T2}
Let $\kappa$ satisfy condition (\ref{As1}).

  \noindent
   $ (i)$ If $ \| T_{\kappa}\|_{HS}<1$ then at least for some $z>1$ there
   is a finite function  $f\geq 1$
    which satisfies
   equation (\ref{A3}).

\noindent
   $ (ii)$ If $ \| T_{\kappa}\|<1$ and  kernel  $\kappa$ satisfies 
condition $(C1)$ or condition $(C3)$ from Theorem \ref{C1},
   then at least for some $z>1$ there
   is a finite function  $f\geq 1$ which satisfies
   equation (\ref{A3}).

\end{theo}

\noindent
{\bf Proof.}
To prove $(i)$ we shall  construct a function $f\geq 1$ which satisfies conditions of
Lemma \ref{L2}.
Let $\| T_{\kappa}\|_{HS} =\lambda < 1$.
Then 
\begin{equation}\label{psi}
  \| \psi\|_2=\lambda
\end{equation}
(see definition of $\psi$ in (\ref{dpsi})).
For any $\varepsilon\geq 0$ let us define
\begin{equation}\label{g}
g(x,\varepsilon)=  T_{\kappa}\left[ e^{\varepsilon \psi}-1\right](x)
=\int_{S} \kappa  (x,y)\left(e^{\varepsilon \psi(y)}-1\right)
d \mu(y).
\end{equation}
By the Cauchy-Bunyakovskii inequality 
\begin{equation}\label{g1}
g(x,\varepsilon)\leq  
\left(\int_{S} \kappa ^2 (x,y)
d \mu(y)\right)^{1/2}
\left(\int_{S} \left(e^{\varepsilon \psi(y)}-1\right)^2
d \mu(y)\right)^{1/2}= \psi(x)A(\varepsilon),
\end{equation}
where function
\[A(\varepsilon):=
\left(\int_{S} \left(e^{\varepsilon
      \psi(y)}-1\right)^2 d \mu(y)
\right)^{1/2}\]
 is increasing and by the assumption (\ref{As1}) and the dominated
 convergence  is continuous on
$[0,a/4]$. Furthermore,  for $0<\varepsilon<a/4$ we can compute
 \[A'(\varepsilon)=
 \frac{\int_{S} \psi(y)e^{\varepsilon
      \psi(y)}\left(e^{\varepsilon
      \psi(y)}-1\right)d \mu(y)
}{
  \left( \int_{S} \left(e^{\varepsilon \psi(y)}-1\right)^2 d \mu(y)\right)^{1/2}
}.\]
Using again the Cauchy-Bunyakovskii inequality and condition (\ref{As1}) we derive from here
that for any small positive $\varepsilon$
 \[A'(\varepsilon)\leq 
  \left( \int_{S}\psi^2(y)e^{2\varepsilon
      \psi(y)}d \mu(y)\right)^{1/2}
\leq 
  \left( \int_{S}Me^{3\varepsilon
      \psi(y)}d \mu(y)\right)^{1/2}<\infty,\]
  where $M$ is some absolute positive constant.
  Hence, taking into account 
   (\ref{psi}) we have
   \[\limsup_{\varepsilon \downarrow 0}A'(\varepsilon) \leq  \|
   \psi\|_2=\lambda< \lambda +\frac{1-\lambda}{2}
 =:  \lambda_1<1.\]
This bound together with  $A(0)=0$  and the mean-value Theorem allows us to conclude that
there exists some
positive value $\varepsilon_0>0$ such that for all $0<\varepsilon<\varepsilon_0$
\begin{equation}\label{g3}
A(\varepsilon)<\lambda_1 \varepsilon.
\end{equation}
Therefore for all $0<\varepsilon<\varepsilon_0$ we get by (\ref{g1})
\begin{equation}\label{A4}
g(x,\varepsilon)\leq  
\psi(x)A(\varepsilon)<\lambda_1 \varepsilon \psi(x).
\end{equation}

 Now fix $z>1$ arbitrarily and denote
 ${\widetilde \psi}= z \psi$. 
Define also a function 
\begin{equation}\label{g1'}
{\widetilde g}(x,\varepsilon)=
z T_{\kappa}\left[ e^{\varepsilon {\widetilde \psi}}-1\right](x)
= z T_{\kappa}\left[ e^{\varepsilon z \psi}-1\right](x)
=zg(x,z \varepsilon).
\end{equation}
According to 
(\ref{A4}) we have
\begin{equation}\label{A5}
{\widetilde g}(x,\varepsilon) \leq z\varepsilon \lambda_1 {\widetilde
  \psi}(x)
\end{equation}
for all  $0< \varepsilon<\varepsilon_0$.
Let us set
\[f_z=z\left( e^{\varepsilon {\widetilde
  \psi}}-1\right)+1.\]
We claim, that for some $z>1$
\begin{equation}\label{A15}
  \Phi_z[f_z]\leq f_z.
\end{equation}
Indeed,  consider
\begin{equation}\label{A6}
\Phi_z[f_z]:=  \Phi_{z, \kappa}[f_z]=ze^{T_{\kappa}[f_z-1]}=ze^{zT_{\kappa}[
     e^{\varepsilon {\widetilde
  \psi}}-1]}.
\end{equation}
Using  definition (\ref{g1'}) and bound (\ref{A5}) we obtain
from here
\begin{equation}\label{A7}
  \Phi_z[f_z](x)=ze^{{\widetilde g}(x,\varepsilon)} \leq ze^{z\varepsilon \lambda_1 {\widetilde
  \psi}(x)}.
\end{equation}
Let us assume now that $1<z<\delta /\lambda_1$ for some
$\lambda_1<\delta <1$.
Then we have
\begin{equation}\label{A8}
e^{z\varepsilon \lambda_1 {\widetilde
  \psi}(x)} \leq e^{\varepsilon \delta  {\widetilde
  \psi}(x)}.
\end{equation}
Under  assumption (\ref{inf}) we  have
$\psi(x)>b>0$ for some positive $b$, which implies that
${\widetilde
  \psi}(x)>b$ as well. Therefore
we can find  $1<z<\delta /\lambda_1$
such that for all $x\in S$
\[
e^{\varepsilon \delta  {\widetilde
  \psi}(x)}\leq e^{\varepsilon {\widetilde
  \psi}(x)}-\frac{z-1}{z},
\]
which  together with (\ref{A8}) gives us
\[ze^{z\varepsilon \lambda_1 {\widetilde
  \psi}(x)} \leq
ze^{\varepsilon \delta {\widetilde
  \psi}(x)}\leq z
\left(  e^{\varepsilon {\widetilde
  \psi}(x)}-\frac{z-1}{z}\right)=
z\left( e^{\varepsilon {\widetilde
  \psi}(x)}-1\right)+1=f_z(x).
\]
Substituting this bound into (\ref{A7}) we finally get (\ref{A15}). 
Hence, function $f_z$ satisfies the conditions of Lemma \ref{L2}, by
which the statement $(i)$ of
Theorem \ref{T2} follows.

The proof of statement $(ii)$
is  very similar to the previous one.
Let $\|T_{\kappa}\|=\lambda<1$. Assume first that condition $(C3)$
is satisfied.
Recall that by
Lemma 5.15 \cite{BJR} 
operator $T_{\kappa}$ with a finite Hilbert-Schmidt norm (assumption (\ref{HS})) has a nonnegative
eigenfunction $\phi \in
L^2(S, \mu)$ such that $T_{\kappa}\phi=\|T_{\kappa}\|\phi.$ Hence,
there is a function $\phi $ such  that $\|\phi\|_2=1$ and
\begin{equation}\label{I}
  \phi(x)=\frac{1}{\lambda}\int_{S} \kappa  (x,y)
   \phi(y)
d \mu(y).
\end{equation}
This together with the Cauchy-Bunyakovskii inequality
immediately implies
\begin{equation}\label{I1}
  \phi(x)\leq \frac{1}{\lambda} \psi(x).
\end{equation}
Under the assumptions of monotonicity of $\kappa$ equality
(\ref{I})
implies that $\phi(x)$ is also monotone increasing, therefore
\begin{equation}\label{S1}
  \phi(x)=\frac{1}{\lambda}\int_{S} \kappa  (x,y)
   \phi(y) d \mu(y)
   \geq \frac{1}{\lambda} T_{\kappa}[1] (x) \int_{S} \phi(y)d \mu(y)
= \frac{c_2}{\lambda} T_{\kappa}[1] (x) 
\end{equation}
for some $c_2>0$.
Next, taking into account condition (\ref{TT}) we  derive
\begin{equation}\label{S2}
\psi^2(x) =\int_{S} \kappa  ^2(x,y)
   d \mu(y)\leq c_1^2 (T_{\kappa}[1] (x))^2
\int_{S} (T_{\kappa}[1] (y))^2
   d \mu(y)\leq c_1^2 (T_{\kappa}[1] (x))^2 \|T_{\kappa}\|_{HS}^2.
   \end{equation} 
Combining now  (\ref{S2}),  (\ref{S1}) and  (\ref{I1}) we get
\[
 \frac{c_2}{\lambda} T_{\kappa}[1] (x) \leq  
 \phi(x)\leq \frac{1}{\lambda} \psi(x)
\leq \frac{1}{\lambda} c_1 T_{\kappa}[1] (x) \|T_{\kappa}\|_{HS},
  \]
  which immediately yields
\begin{equation}\label{S3}
 \frac{ \psi(x) }{\phi(x) }\leq \frac{ c_1 \lambda\|T_{\kappa}\|_{HS} }{c_2 }
   \end{equation}
   for all $x\in S$.
Notice that (\ref{inf})  implies
\begin{equation}\label{I3}
  \phi(x)\geq  c_0 >0
\end{equation}
for all $x\in S$ and some $c_0>0$.

We can  show now that
function
\begin{equation}\label{ReF}
  F_z(x)=z\left(e^{\varepsilon \phi(x)}-1\right)-1
\end{equation}
satisfies the conditions of Lemma \ref{L2} for some positive $\varepsilon$. First we consider similar
to (\ref{g}) function
\[
G(x,\varepsilon)=  T_{\kappa}\left[ e^{\varepsilon \phi}-1\right](x)
=\int_{S} \kappa  (x,y)\left(e^{\varepsilon \phi(y)}-1\right)
d \mu(y).
\]
From here we derive using assumption (\ref{As1}) and bound (\ref{I1}), 
that  at least for all $\varepsilon<a \lambda /4$
functions 
\begin{equation}\label{g'}
\frac{\partial}{\partial \varepsilon}G(x,\varepsilon)
=\int_{S} \kappa  (x,y) \phi(y)e^{\varepsilon \phi(y)}
d \mu(y) 
\end{equation}
and
\begin{equation}\label{g''}
 \frac{\partial^2}{\partial \varepsilon^2}G(x,\varepsilon)
= \int_{S} \kappa  (x,y) \phi(y)^2 e^{\varepsilon \phi(y)}
d \mu(y) 
\end{equation}
are  finite and non-negative for any $x\in S$.
Note that for all $x\in S$
\[G(x,0)=0,\]
and
\begin{equation}\label{g'0}
\frac{\partial}{\partial \varepsilon}G(x,\varepsilon)\mid _{\varepsilon=0}
=\int_{S} \kappa  (x,y) \phi(y)
d \mu(y) =\lambda \phi(x).
\end{equation}
Therefore for all $x\in S$ and $0\leq
\varepsilon<a \lambda /4$ we have
\begin{equation}\label{S4}
G(x,\varepsilon) \leq \varepsilon \left( \lambda \phi(x)
+\varepsilon \frac{\partial^2}{\partial \varepsilon^2}G(x,\varepsilon)
\right).
\end{equation}
Under the assumption (\ref{As1}) we get from (\ref{g''}) and
(\ref{I1}),
that for all
$0\leq
\varepsilon<a \lambda /4$
\begin{equation}\label{g''1}
 \frac{\partial^2}{\partial \varepsilon^2}G(x,\varepsilon)
\leq \psi(x) \int_{S}  \phi(y)^4 e^{2\varepsilon \phi(y)}
d \mu(y) \leq  \psi(x) c_3,
\end{equation}
where $c_3$
is some positive constant. Taking also  into account bound (\ref{S3}),
we derive from (\ref{g''1})
\[\frac{\partial^2}{\partial \varepsilon^2}G(x,\varepsilon)
\leq c_4 \phi(x)\]
for some positive constant $c_4$.
Substituting this into (\ref{S4}), we get
\begin{equation}\label{A12}
G(x,\varepsilon) \leq \varepsilon \phi(x) \left( \lambda 
+\varepsilon c_4 \right).
\end{equation}
It is clear that
for all
 small $\varepsilon >0$ we have
\[\lambda 
+\varepsilon c_4 <\lambda +\frac{1-\lambda}{2}:=\lambda
_1<1.\]
This together with (\ref{A12}) immediately yields
\begin{equation}\label{A14}
G(x,\varepsilon) \leq \lambda _1 \varepsilon \phi(x)
\end{equation}
 for all
 small $\varepsilon >0$.
Then repeating almost exactly the same argument which led from (\ref{A4}) to (\ref{A15}),
one can derive from (\ref{A14}) that for $F_z$
defined by (\ref{ReF})
\[
  \Phi_z[F_z]\leq F_z.
\]
Hence, function $F_z$ satisfies the conditions of Lemma \ref{L2},
which yields statement $(ii)$ when $\kappa $ satisfies $(C3)$. 

Finally, let $\|T_{\kappa}\|=\lambda<1$, and let condition $(C1)$
be satisfied. Taking into account assumptions (\ref{inf}) we easily derive bounds similar to (\ref{TT}) and (\ref{S3}). Then the rest of the proof simply repeats the previous one. This completes the proof of statement $(ii)$ and finishes the proof of Theorem \ref{T2}. \hfill$\Box$

\bigskip
\section{Proofs of the main results.}
\setcounter{equation}{0}

\subsection{Proof of Theorem \ref{T1}.} The 
  statement follows immediately by 
 Theorem \ref{T3} and Remark \ref{R1}. \hfill$\Box$

\subsection{Proof of Corollary \ref{C4}.} 

\begin{lem}\label{LRev}
If $\| T_{\kappa}\| = 1$ one has
\begin{equation}\label{N5}
\lim_{c \uparrow 1} r_{c \kappa}=1.
\end{equation}
\end{lem}

\noindent
{\bf Proof.} Let
${\cal X}^c$ denote  the total progeny of the  $B_{c \kappa}$ (see Definition \ref{d1}). 
It is clear that 
if $0<c<c'$ then ${\cal X}^{c'}$ stochastically dominates  ${\cal
  X}^{c}$, and therefore it is obvious 
that $r_{c \kappa}$ is a
monotone 
non-increasing function in $c>0$. Also, it follows from the definition
of  $r_{c \kappa}$, that $r_{c \kappa}\geq 1$. Hence,
there exists $\lim_{c \uparrow 1} r_{c \kappa}\geq 1$.
Assume,
\begin{equation}\label{N6}
\lim_{c \uparrow 1} r_{c \kappa}=:r > 1.
\end{equation}
Define
\[1< z:= \frac{1}{2}+ \frac{r}{2}<\lim_{c \uparrow 1} r_{c \kappa}.\]
Then by Theorem \ref{T1}
for any fixed $c<1$ 
there
exists minimal solution 
  $1 \leq f<\infty$ to (\ref{A3}):
  \[f=ze^{ T_{c\kappa}[f-1]}.\]
Notice that also
  \[f(x)\geq z\]
  for all $x\in S$. Let $c'=\sqrt{z} >1$ and set
\begin{equation}\label{defg}
  g:=\frac{f}{c'}\geq\sqrt{z}> 1.
\end{equation}
It is straightforward to derive
\[\Phi_{\sqrt{z},c'c\kappa}[g]=  \sqrt{z} e^{c'T_{c\kappa}[g-1]}
=  \frac{1}{\sqrt{z}} \, z
e^{T_{c\kappa}[c'g-1]-(c'-1)T_{c\kappa}[1]}= \frac{1}{c'} f 
e^{-(c'-1)T_{c\kappa}[1]}\leq g.\]
Hence, by  Lemma \ref{L2} there exists a  function $1<h<\infty$ such that
\begin{equation}\label{N7}
h= \Phi_{\sqrt{z},c'c\kappa}[h] \equiv \sqrt{z}e^{
  T_{\sqrt{z}c\kappa}[h-1]}.
\end{equation}
Choose now $\frac{1}{\sqrt{z}}<c<1$. Then 
existence of an $a.s.$ finite solution $h>1$ to 
equation (\ref{N7}) with $\sqrt{z}> 1$ implies by 
 Theorem \ref{T1}  that $r_{\sqrt{z}c\kappa}> 1$ even when
$\|T_{\sqrt{z}c\kappa}\| =\sqrt{z}c> 1$. We get a contradiction with
(\ref{rRev}), which finishes the proof of Lemma. \hfill$\Box$

\medskip

By Lemma \ref{LRev}  we have
\[1=\lim_{c \uparrow 1} r_{c \kappa}\geq r_{ \kappa}\geq 1,\]
which yields
\begin{equation}\label{(a)}
r_{\kappa}=1, \ \ \ \mbox{ if } \ \ \ \| T_{\kappa}\| = 1.
\end{equation}
This together with (\ref{rRev}) completes the proof of Corollary \ref{C4}.
  \hfill$\Box$
 
\subsection{Proof of Theorem \ref{C1}.} The 
  statement follows immediately by 
 Theorem \ref{T1} and Theorem \ref{T2}. \hfill$\Box$
 
\subsection{Proof of Proposition \ref{Prev}.}

 Assume, equation (\ref{+})  has an $a.s.$ finite solution
$f$ such that $\sup_{x\in S}f(x)<0$. We shall show that in this case there is $z>1$ such that
equation (\ref{-}) has $a.s.$ finite solution $g\geq 1$. This
by Theorem \ref{T1} will imply  $r_{\kappa} >1$.

By our assumption
\[
f=1-e^{-T_{\kappa}[f]}.\]
Then for $h:=-(f-1)>1$ we have
\begin{equation}\label{Rev1}
h=e^{T_{\kappa}[h-1]}.
\end{equation}
{\it Claim.} There are  $0<\varepsilon<1$ and $z>1$ such that
function
\begin{equation}\label{Revg}
H=\varepsilon +(1-\varepsilon)h
\end{equation}
satisfies inequality
\begin{equation}\label{Rev2}
\Phi_z[H]:=ze^{T_{\kappa}[H-1]}\leq H.
\end{equation}

\noindent
{\it Proof of the Claim.} By (\ref{Rev1}) we have for any $z>1$
\begin{equation}\label{Rev3}
  \Phi_z[H]= ze^{T_{\kappa}[
\varepsilon +(1-\varepsilon)h
  -1]}
= z\left(
  e^{T_{\kappa}[h-1]}\right)^{1-\varepsilon}=zh^{1-\varepsilon}
\equiv z \ \frac{h^{1-\varepsilon}}{\varepsilon +(1-\varepsilon)h}\  H.
\end{equation}
Define  for all $q\geq 1$
\[Q(\varepsilon, q):=\frac{q^{1-\varepsilon}}{\varepsilon
  +(1-\varepsilon)q}.\]
It is straightforward to compute that for any $0<\varepsilon<1$
and for any $q> 1$
\begin{equation}\label{Rev19}
\frac{\partial}{\partial q}Q(\varepsilon, q)=
\frac{(1-\varepsilon)q^{-\varepsilon} (\varepsilon +(1-\varepsilon)q)
  -q^{1-\varepsilon}(1-\varepsilon)}{(\varepsilon
  +(1-\varepsilon)q)^2} < 0.
\end{equation}
Recall that by the assumption $$h_*=\inf_{x\in S}h(x)=1-\sup_{x\in
  S}f(x)>1.$$
Hence, by (\ref{Rev19}) for all $x\in S$
\[
Q(\varepsilon, h(x))\leq Q(\varepsilon, h_*) < Q(\varepsilon, 1)=1.
\]
Setting now
$z=\frac{1}{Q(\varepsilon, h_*)}>1$
we derive from (\ref{Rev3})
\[ \Phi_z[H]\leq H,\]
which proves our Claim.\hfill$\Box$

Notice, that by the definition 
\[
H=\varepsilon +(1-\varepsilon)h =
1+(1-\varepsilon)(h-1)\geq 1,
\]
since $h>1$. This and (\ref{Rev2}) confirm that the conditions on Lemma \ref{L2} are
fulfilled by function $H$. Therefore by Lemma  \ref{L2} there exists
$a.s.$
finite solution $g\geq 1$ to (\ref{-}) with some $z>1$. This completes the proof of
Proposition \ref{Prev}. \hfill$\Box$

\subsection{Proof of Theorem \ref{T1A}.}

\subsubsection{The upper bound.}\label{UB}
\begin{theo}\label{Tub}
If $\|T\|<1$ then under
conditions of Theorem \ref{T1A} one has $r_{\kappa}>1$ and 
\begin{equation}\label{lb}
\lim_{n \rightarrow \infty}{\bf P}\left\{  C_1 \Big(
G^{\cal V}(n,\kappa )
    \Big) > \left( \frac{1}{\log r_{\kappa}} + \delta \right)\log n \right\}
= 0.
\end{equation}
for any $\delta >0$.
\end{theo}

\noindent{\bf Proof.} Notice that here $r_{\kappa}>1$  simply by 
Theorem
\ref{C1}.

Throughout the proof we assume that the condition $(C3)$ is satisfied. In
 the case of $(C1)$ the proof repeats the same arguments with obvious simplifications. 

Recall the usual algorithm of finding a
connected component in a random graph. Given the sequence 
$(x_1, \ldots , x_n)$ and a corresponding graph
$ G^{\cal V}(n,\kappa)$, take any vertex $1\leq i \leq n$ to be the root. Find all
vertices  connected to this vertex $i$ in the graph $ G^{\cal V}(n,\kappa)$, and
then mark $i$ as "saturated". Then for each
non-saturated  revealed vertex, we find all  vertices
connected to it but which have not been used previously. 
We continue this process until we end
up with a tree of saturated vertices.

Denote  $\tau_{n}^i$ the set of vertices in
the tree constructed according to the  above algorithm with the
root at  vertex $i$.
Then for any $\omega>0$
\begin{equation}\label{ct}
{\bf P}\left\{  C_1 \Big(
G^{\cal V}(n,\kappa )\Big)
>\omega\right\} =
{\bf P} \left\{ \max _{1\leq i\leq n} |\tau_{n}^i| > \omega
\right\} .
\end{equation}

Let constant $a$
be the one from condition (\ref{As1}). Then for any
\begin{equation}\label{q}
0\leq q<a/2
\end{equation}
define an auxiliary
probability measure on $S$:
\begin{equation}\label{muq}
  \mu_q(x) = m_q e^{qT_{\kappa}[1](x)}\mu(x)
\end{equation}
with a normalizing constant
\[m_q:=\left(\sum_S e^{qT_{\kappa}[1](x)}\mu(x) \right)^{-1},\]
which is strictly positive due to assumption (\ref{As1}).
Notice that $\mu_0(x) =\mu(x)$ for all $x\in S$, and $m_q$ is continuous in $q$ on $[0, a/2]$ with 
$m_0=1$. This implies in particular, that for any $\varepsilon '>0$ one can choose positive  $q$ so that
\begin{equation}\label{Aex}
\mu(x)\leq (1+\varepsilon' )\mu_q (x)
\end{equation}
for all $x$.
Fix $\varepsilon>0$ and $0<q<a/2$ arbitrarily and define an event
\begin{equation}\label{calB}
{\cal B}_n =
\left\{   \frac{\#\{1\leq i\leq n: x_i=x\}}{n}-\mu(x) \leq
  \varepsilon \mu_q(x), \ \ \mbox{ for all } \ \  x\in S \right\}.
\end{equation}
By 
 assumption (\ref{As3})
 we have
\begin{equation}\label{Re11}
{\bf P} \left\{{\cal B}_n \right\}=1-o(1).
\end{equation}
Then we  derive from  (\ref{ct})  that
\begin{equation}\label{Rev23}
 {\bf P} \left\{  C_1\Big(  G^{\cal V}(n,\kappa )   \Big) > \omega
\right\}
\leq
{\bf P} \left\{ \max _{1\leq i\leq n} |\tau_{n}^i| >  \omega \mid {\cal B}_n
\right\} +o(1).
\end{equation}
Notice that  the distribution of the size
$|\tau_{n}^i|$ depends only on
the type $x_i$ of vertex $i$. Then using  notation
\begin{equation}\label{Ren}
|\tau_{n}(x)|=_d |\tau_{n}^i| \Big| _{x_i=x}
\end{equation}
for each $x\in S$, we derive from (\ref{Rev23})
\begin{equation}\label{SA18}
 {\bf P} \left\{  C_1\Big(  G^{\cal V}(n,\kappa )   \Big) > \omega
\right\}
\leq n\sum_{x\in S}
\left({\mu}(x)+\varepsilon {\mu}_q(x)\right)  
{\bf P} \left\{  
|\tau_{n}(x)| >  \omega\mid {\cal B}_n 
\right\} +o(1)
\end{equation}
as $n \rightarrow \infty$.

To approximate the distribution of $|\tau_{n}(x)|$
we shall use the following branching processes.
For any  $c\geq 1$ and $ q\geq 0$  let
$B_{c, q}$ denote a process defined similar to
$B_{\kappa}$ in Definition \ref{d1},
but with
the distribution of the offspring 
\[Po\left(c\kappa(x,y)\mu_q(y)\right)\]
instead of $Po\left(\kappa(x,y)\mu(y)\right)$. Notice, that
$B_{1, 0}$
is  defined exactly as  $B_{\kappa}$.
Let further ${\cal X}^{c,q}(x)$ denote the total number of the particles
(including the initial one) produced by
the branching process $B_{c, q}$
starting with a single particle
of  type $x$.

\begin{prop}\label{LJ1} Under conditions of Theorem \ref{T1A} one can find 
$q>0$ and
  $c>1$ arbitrarily close to $0$ and $1$,
  correspondingly, such that for some $\varepsilon>0$ in the
  definition of ${\cal B}_n$ 
\begin{equation}\label{SA23}
{\bf P} \left\{ |\tau_n(x) | > \omega \mid {\cal B}_n
\right\} 
\leq
{\bf P} \left\{
  {\cal X}^{c,q} (x)> \omega
\right\} 
\end{equation}
for all $x\in S$, $\omega >0$, and for all large $n$.
\end{prop}

\noindent
{\bf Proof.}
Observe that at each step of the exploration algorithm which defines $\tau_n^i$,
the number of the type $y$ offspring 
of a particle of type $x$ has a
binomial distribution $ Bin(N_y',{ p}_{xy}(n))$ where $N_y'$
is the number of the remaining vertices of type $y$.

We shall use a well-known fact that 
 a binomial $Bin(n,p)$ distribution is dominated by
 a Poisson  distribution $Po(-n\log(1-p))$. 
First we shall derive an upper bound for $N_y'$.
Notice that conditionally on ${\cal B}_n$ we have
\begin{equation}\label{Ny}
N_y'\leq \#\{1\leq i\leq n: x_i=y\}\leq n(\mu(y)+\varepsilon \mu_q(y))
\end{equation}
for each $y\in S$. The last inequality implies that for any $y$ such that
\[\#\{1\leq i\leq n: x_i=y\}>0\]
we have
\begin{equation}\label{Fe1}
 n( \mu(y)+\varepsilon \mu_q(y)) \geq 1.
\end{equation}
By the Cauchy-Bunyakovskii inequality and by assumption  (\ref{As1}) we have
\[
\sum_S e^{qT_{\kappa}[1](x)}\mu(x)\leq \sum_S e^{q\psi(x)}\mu(x)<\infty
\]
for all $q\leq a$. Hence,  for all $0<q<a/2$ and for all $y\in S$
\[
 \mu(y)+\varepsilon \mu_q(y) \leq
b_3\left(e^{-a T_{\kappa}[1](y)}+\varepsilon m_q e^{(q-a)T_{\kappa}[1](y)} \right)\leq
 b_2e^{-a T_{\kappa}[1](y)/2},
\]
for all  $0<\varepsilon<1$,
where  $b_2, b_3$ are some positive constants. 
Combining this with (\ref{Fe1}) we obtain for  all
$y$ such that
$\#\{1\leq i\leq n: x_i=y\}>0$ 
\[\frac{1}{n}\leq \mu(y)+\varepsilon \mu_q(y) \leq b_2e^{-aT_{\kappa}[1](y)/2}.\]
This implies that conditionally on ${\cal B}_n$ 
\[
\max_{x \in \{x_1, \ldots, x_n\}} T_{\kappa}[1](x)\leq A_1 \log n
\]
for some constant $A_1.$
Taking into account assumption (\ref{TT}), we derive from here that
for all large $n$ conditionally on ${\cal B}_n$ 
\begin{equation}\label{Kn1}
{ p}_{x_i x_j}(n) = \frac{\kappa(x_i, x_j)}{n}\leq \frac{c_1(A_1\log n)^2}{n}.
\end{equation}
The last bound and (\ref{Ny}) together with (\ref{Aex}) allow us 
for any fixed 
positive $\varepsilon _1$
to choose $\varepsilon$ and $q$  so that
conditionally on ${\cal B}_n$ we get
\[
    -N_y' \log{(1-{ p}_{xy}(n))}
 \  \leq \left({\mu}(y)+\varepsilon {\mu}_q(y)\right)
 n |\log{(1-{ p}_{xy}(n))}|
\]
\begin{equation}\label{pB1}
 \leq (1+\varepsilon _1){\mu}_{q}  (y)\kappa(x,y)=:c \mu_{q}(y)\kappa(x,y)
\end{equation}
for all large $n$ and all $x,y \in \{x_1, \ldots, x_n\}$. Hence, (\ref{pB1}) holds for any
$q>0$ and $c>1$ arbitrarily close to $0$ and $1$, respectively.
It follows by (\ref{pB1}) that the
 binomial distribution
$Bin(N_y' ,{ p}_{xy}(n))$ is dominated stochastically by the Poisson distribution $Po(c\mu_{q}  (y)\kappa (x,y))$.
Therefore if conditionally on ${\cal B}_n $
 at each  step
of the exploration algorithm which reveals $\tau_n^i$, we replace the
$ Bin(N_y',{ p}_{xy}(n))$ variable 
with the
$Po\left(c\mu_{q}  (y) \kappa (x,y)\right)$ one, we arrive at
the statement (\ref{SA23}) of the
 Proposition. \hfill$\Box$

\bigskip
Substituting (\ref{SA23})  into (\ref{SA18})
we derive that for any 
$q>0$ and $c>1$, one has
\[
{\bf P} \left\{   C_1\Big(  G^{\cal V}(n,\kappa )  \Big) >\omega \right\} 
\leq
bn
\sum_{x\in S} 
 \mu _{q}(x)
{\bf P} \left\{ {\cal X}^{c,q} (x)> \omega
\right\}
 +o(1)
\]
as $n \rightarrow \infty$, where
$b$ is some positive constant. 
This bound together with  the Markov inequality  imply for all
$z\geq 1$
\begin{equation}\label{SA25}
{\bf P} \left\{   C_1\Big(  G^{\cal V}(n,\kappa )  \Big) >\omega \right\} 
\leq b n z^{-\omega}
\sum_{x\in S} 
 \mu _{q}(x){\bf
    E} z^{{\cal X}^{c,q} (x)}
+o(1).
\end{equation}

Let $T_{c\kappa, \mu_q}$ denote an integral operator associated with
branching process $B_{c,q}$
\begin{equation}\label{defT}
T_{c\kappa, \mu_q}[f](x)=\int_S c\kappa(x,y)f(y)d\mu_q(y)=
\sum_S c\kappa(x,y)f(y)\mu_q(y).
\end{equation}
Assume from now on that $q>0$ and $c>1$  are such that
\begin{equation}\label{S21}
cm_q \geq 1.
 \end{equation}
We shall extend now the result from Lemma 7.2
in \cite{BJR} on the approximation of kernels
for our special case of  unbounded kernels. 
First,  taking into account conditions  (\ref{As1}) and  (\ref{TT})
we derive that for any fixed $q<a/4$ and $c>1$
\[\|T_{c\kappa, \mu_q}\|_{HS}^2  = \int_S \int_S (cm_q)^2
\kappa^2(x,y)e^{q T_{\kappa}[1](x)}e^{q T_{\kappa}[1](y)}
d\mu(x)d\mu(y)<\infty.\]
Then by the Theorem on
Dominated Convergence
\begin{equation}\label{S22}
 \|T_{c\kappa, \mu_q}-T_{\kappa}\|_{HS}^2  \rightarrow 0
 \end{equation}
 as $c\rightarrow 1$ and $q \rightarrow 0$.
Furthermore, since
\[\|T_{\kappa}\|\leq  \|T_{c\kappa, \mu_q}\| \leq  \|T_{\kappa}\|
+\|T_{c\kappa, \mu_q}-T_{\kappa}\|_{HS},
\]
convergence (\ref{S22}) implies as well $$\|T_{c\kappa,
  \mu_q}\|\rightarrow
\|T_{\kappa}\|
$$
as $c\rightarrow 1$ and $q \rightarrow 0$.
Hence, if $\|T_{\kappa}\|<1$ then we can choose  $0<q<a/4$ and $c>1$
so that (\ref{S21}) holds  together with
\begin{equation}\label{A22}
 \|T_{c\kappa, \mu_q}\|<1.
\end{equation}
Now
 for all values $c$ and $q$ for which (\ref{A22}) holds
 we have by Theorem  \ref{T2},
part $(ii)$ that
\begin{equation}\label{N}
r(q,c):=\sup \{z\geq 1: \sum_{x\in S}{\bf E}z^{{\cal X}^{c,q} (x)} \
\mu_q(x) <\infty\} >1,
\end{equation}
and therefore 
for all $1<z < r({q,c})$ 
\begin{equation}\label{Ma2}
\sum_{x\in S} 
 \mu _{q}(x)
 {\bf E} z^{{\cal X}^{c,q} (x)} <\infty .
\end{equation}

Notice, that condition (\ref{S21})  implies that
${\cal X}^{c,q}(x)$ is stochastically larger than ${\cal X}(x)$ for
any $x\in S$, which clearly yields 
\begin{equation}\label{S27}
r(q,c)\leq r_{\kappa}.
\end{equation}

\begin{lem}\label{LS1}
For any $z<r_{\kappa}$ there are $q>0$ and $c\geq 1/m_q$ such that
\begin{equation}\label{S23}
z\leq r(q,c)\leq r_{\kappa}.
\end{equation}
  \end{lem}

  \noindent{\bf Proof.}
Notice that when $z\leq
1$ statement (\ref{S23}) follows  by  (\ref{N}) and (\ref{S27}).

Let us fix  $1<z<r_{\kappa}$ arbitrarily.
We shall show that for some $q>0$ and  $c\geq 1/m_q$ 
equation
\begin{equation}\label{RevS24}
f= ze^{T_{c\kappa, \mu_q}[f-1]} 
\end{equation}
has 
a finite solution $f\geq 1$. This by Theorem \ref{T1} will imply that
$
z \leq r(q,c).
$
The later together with (\ref{S27}) would  immediately yield (\ref{S23}). 

First we rewrite equation (\ref{RevS24}).
Let $q>0$ and $c\geq 1/m_q$  be such
that (\ref{A22}) holds. Set
$$c_q:=c m_q \geq
1,$$
and denote 
\[\widetilde{\kappa}(x,y)= c_q{\kappa}(x,y)e^{qT_{\kappa}[1](y)}.\]
Then  (\ref{RevS24}) becomes
\begin{equation}\label{S24}
f=  \Phi_{z,\widetilde{\kappa}}[f].
\end{equation}
Since
\[\Phi_{z,\widetilde{\kappa}}[f]=z\exp\{T_{\kappa}[c_qe^{qT_{\kappa}[1]}f]-
T_{\widetilde{\kappa}}[1]\},
\]
setting $g=c_qe^{qT_{\kappa}[1]}f$ we can rewrite equation (\ref{S24})  as follows
\[
g=c_q z\exp\{T_{\kappa}[g-1] + T_{\kappa}[1] + q  T_{\kappa}[1] 
-
T_{\widetilde{\kappa}}[1]\}
\]
\[= \Phi_{c_q z,\kappa}[g] e^{(1+q)T_{\kappa}[1]-
T_{\widetilde{\kappa}}[1]}.\]
Hence, equation (\ref{S24}) has 
a finite solution $f\geq 1$ if and only if equation
\begin{equation}\label{S30}
g=\Phi_{c_q z,\kappa}[g] e^{(1+q)T_{\kappa}[1]-
T_{\widetilde{\kappa}}[1]}=:G[g]
\end{equation}
has 
a finite solution $g\geq c_qe^{qT_{\kappa}[1]}$.
Observe that $G$ is  a monotone operator, i.e., if $g\geq g_1$ then 
$G[g]\geq G[g_1]$. Since
\[G[g]= c_qe^{qT_{\kappa}[1]}\Phi_{z,\widetilde{\kappa}}[c_q^{-1}
e^{-qT_{\kappa}[1]}g], \]
for any
\[
g\geq c_qe^{qT_{\kappa}[1]}
\]
we have
\[G[g]\geq c_qe^{qT_{\kappa}[1]}.\]
If we find a function $g_0$ such that
\begin{equation}\label{O1}
g_0\geq c_qe^{qT_{\kappa}[1]}
\end{equation}
and
\begin{equation}\label{O2}
G[g_0]\leq g_0,
\end{equation}
we can  derive
(using the argument similar to the proof of Lemma 2.1) 
that 
\begin{equation}\label{S}
g:=
\lim_{n\rightarrow \infty}
G^n[g_0]\geq c_qe^{qT_{\kappa}[1]}
\end{equation}
  is the  finite solution to (\ref{S30}).

  Let $g_0$ be the minimal positive solution to 
\begin{equation}\label{S29}
g_0= \Phi_{c_qz,\kappa}[g_0],
\end{equation}
where we assume that
\begin{equation}\label{O3}
c_qz < r_{\kappa}.
\end{equation}
By Theorem \ref{T1}
the minimal positive solution to (\ref{S29}) is finite. Furthermore, according
to formula (\ref{O10}) we have
\begin{equation}\label{O4}
g_0\geq \Phi_{c_qz,\kappa}^2[1]=c_qz e^{(c_qz
  -1)T_{\kappa}[1]}\geq
c_q e^{(c_q
  -1)T_{\kappa}[1]},
\end{equation}
where we used the fact that $z>1$. Now for a fixed previously $1<z<
r_{\kappa}$ we can choose 
\begin{equation}\label{O6}
0<q\leq \frac{r_{\kappa}}{z} -1,
\end{equation}
and then set
\begin{equation}\label{O6a}
c_q=1+q.
\end{equation}
With this choice of constants we have
condition (\ref{O3}) satisfied, and moreover
from (\ref{O4}) we derive
\begin{equation}\label{O5}
g_0\geq c_q e^{(c_q
  -1)T_{\kappa}[1]}=c_q e^{q T_{\kappa}[1]},
\end{equation}
which means that condition  (\ref{O1}) is satisfied as well.
Notice also that by (\ref{O6}) and (\ref{O6a})
\begin{equation}\label{O7}
(1+q)T_{\kappa}[1](x)-
T_{\widetilde{\kappa}}[1](x)=
\int_S\left( 1+q - c_q e^{q T_{\kappa}[1](x)}
\right)\kappa(x,y)d\mu(y)
\end{equation}
\[ \leq
\int_S\left( 1+q - c_q 
\right)\kappa(x,y)d\mu(y) =0.
\]
Therefore   with constants  (\ref{O6}) and (\ref{O6a}) we derive from 
(\ref{S30}), (\ref{O7}) and (\ref{S29}) that
\[
G[g_0]=\Phi_{c_q z,\kappa}[g_0] e^{(1+q)T_{\kappa}[1]-
T_{\widetilde{\kappa}}[1]}\leq \Phi_{c_q z,\kappa}[g_0] =g_0.
\]
Hence,
conditions (\ref{O1}) and (\ref{O2}) are fulfilled. Then  by 
(\ref{S}) equation (\ref{S30}) has a desired finite solution. In turn,
this implies that
equation (\ref{S24}) has a finite solution $f\geq 1$, which yields
statement (\ref{S23}). \hfill$\Box$
\bigskip

By Lemma \ref{LS1} for any $\delta >0$ 
we can choose a small $\delta' >0$ and $(q,c)$
close to $(0,1)$ so that (\ref{Ma2}) holds with 
\[
z=r(q,c)-\delta' > 1,
\]
and moreover
\begin{equation}\label{Ma4}
\left(\frac{1}{\log r_{\kappa}}+\delta \right)\log (r(q,c)-\delta') > 1.
\end{equation}
Now setting
$\omega= \left(\frac{1}{\log r_{\kappa}}+\delta \right)\log n$
and $z=r(q,c)-\delta'$
in (\ref{SA25}) we derive with a help of (\ref{Ma2})
\[
{\bf P} \left\{  C_1 \Big(G^{\cal V}(n,\kappa ) \Big) > \left(\frac{1}{\log r_{\kappa}}+\delta \right)\log n \right\} 
\leq b_1nz^{-\omega}+o(1)
\]
\[
=
b_1 n \exp\left\{-\log (r(q,c)-\delta')\left(\frac{1}{\log r_{\kappa}}+\delta \right)\log n\right\}
 + o(1) \]
where $b_1$ is some finite positive constant. This together with (\ref{Ma4})
yields  statement (\ref{lb}).

\hfill$\Box$

\subsubsection{The lower bound.}
\begin{theo}\label{Tlb}
If $\|T\|<1$ then under 
 conditions of Theorem \ref{T1A}
one has $r_{\kappa}>1$ and 
\begin{equation}\label{ub}
\lim _{n \rightarrow \infty}{\bf P}\left\{  C_1 \Big(
G^{\cal V}(n,\kappa )
    \Big) < \left( \frac{1}{\log r_{\kappa}} - \delta \right)\log n \right\}
= 0
\end{equation}
for any $\delta >0$.
\end{theo}

\noindent
{\bf Proof.}
Fix any small positive $\delta $ and  denote 
\begin{equation}\label{om}
\omega= \left(\frac{1}{\log r_{\kappa}}-\delta\right)\log n, 
\end{equation}
\begin{equation}\label{fN}
N=N(n)= \frac{n}{\omega^2}.
\end{equation}
Introduce also for an arbitrarily fixed finite $D\in S$ and $\varepsilon_1>0$ an
event
\[{\cal A}_n=
\left\{\frac{\#\{x_i :
    x_i=y\}}{n}-\mu(y)\geq - \varepsilon_1 \mu(y), \ \ \mbox{ for all } \
  \ 0 \leq y\leq D \right\} \cap {\cal B}_n
\]
with ${\cal B}_n$ defined by (\ref{calB}).
Observe that by  assumption (\ref{set}) and by (\ref{Re11})
\begin{equation}\label{Rev22}
{\bf P}\left\{ {\cal A}_n
 \right \} \rightarrow 1
\end{equation}
as $n\rightarrow \infty$.
Let
\[{\bf P}_{{\cal A}_n}(\cdot)= {\bf P} \left\{\cdot \mid {\cal A}_n \right\}\]
denote the conditional probability. 

Given a  graph $G^{\cal V}(n,\kappa )$
we shall  reveal recursively  its connected components in the following way.
Let $V_1$ be a random vertex uniformly 
distributed on $\{1, \ldots, n\}$, and let  $L_1=\tau_n^{V_1}$  be
the set of vertices in the tree 
with a root at vertex $V_1$ (see definition of the algorithm in
Section \ref{UB}).

For any $U \subset \{1, \ldots, n\}$ let $\tau_n^{i, U}$
denote a set of vertices of the tree constructed in the same way as 
$\tau_n^i$  but on the set of vertices $\{1, \ldots, n\}\setminus U$ instead of
$\{1, \ldots, n\}$. In particular, with this notation
$\tau_n^{i, \emptyset}=\tau_n^i$.

Given  constructed components $L_1, \ldots , L_k$ for  $1\leq k \leq
[N]$, let 
$V_{k+1}$ be a vertex uniformly distributed on
$\{1, \ldots, n\}\setminus \cup_{i=1}^kL_i$, and set 
$L_{k+1}=\tau_n^{V_k, \cup_{i=1}^kL_i }(V_{k+1})$. If
$\{1, \ldots, n\}\setminus \cup_{i=1}^kL_i=\emptyset$, 
we simply set $L_{k+1}=\emptyset$.
Then according to   (\ref{Rev22})  we have
\begin{equation}\label{F1}
{\bf P}\left\{  C_1 \Big(
G^{\cal V}(n,\kappa )
    \Big) < \left( \frac{1}{\log r_{\kappa}} - \delta \right)\log n \right\}\leq
{\bf P}_{{\cal A}_n} \left\{\max_{1\leq i\leq [N]+1}|L_i|<\omega \right\}
+o(1)
\end{equation}
as $n \rightarrow \infty$.
Consider now
\begin{equation}\label{Re1}
{\bf P}_{{\cal A}_n} \left\{\max_{1\leq i\leq [N]+1}|L_i|<\omega \right\}
\end{equation}
\[={\bf P}_{{\cal A}_n} \left\{|L_{1}|<\omega \right\} 
\prod_{i=1}^{[N]}{\bf P}_{{\cal A}_n} \left\{|L_{i+1}|<\omega \mid
|L_1|<\omega, \ldots , |L_i|<\omega
\right\}.
\]

Similar to (\ref{Ren}), let us  also define  
\[
|\tau_{n}^U(x)|:=_d |\tau_{n}^{i, U}| \Big| _{x_i=x}
\]
for each $x\in S$, $i\in \{1, \ldots, n\}$ and $U\subset \{1, \ldots, n\}$.
Notice  that if $U \subset U' $ then  $|\tau_n^{U'}(x)|$ is
stochastically
dominated by $|\tau_n^U(x)|$ for any $x\in S$. This allows us to
derive from (\ref{Re1}) that
\begin{equation}\label{F4}
{\bf P}_{{\cal A}_n} \left\{\max_{1\leq i\leq [N]+1}|L_i|<\omega \right\}
\leq \left(\sup_{x\in S}  \ \max_{U\subset \{1, \ldots, n\}: |U|\leq  N\omega}
{\bf P}_{{\cal A}_n} \left\{|\tau_n^U(x)|<\omega \right\}
\right)^{N}.
\end{equation}

To approximate the distribution of $|\tau_n^{U}(x)|$ we
introduce another branching process which will be stochastically
dominated by $B_{\kappa}$. 
First define for any value $D\in S$ a probability measure $\hat{\mu}_D$
\begin{equation}\label{muD}
\hat{\mu}_D(y)=
\left\{
\begin{array}{ll}
M_D^{-1}\mu(y), & \mbox{ if } y\leq D, \\
0, & \mbox{ otherwise },
\end{array}
\right.
\end{equation}
where $M_D:=\sum_{y\leq D}\mu(y) $ is a normalizing constant. 
Then for any positive $c$
and $D$ let $\hat{B}_{c, D}$ be a process defined similar to
$B_{\kappa}$,
but with
the distribution of  offspring 
\begin{equation}\label{BD}
Po\left(c\kappa(x,y)\hat{\mu}_D(y)\right)
\end{equation}
instead of $Po\left(\kappa(x,y)\mu(y)\right)$. Notice, that
$\hat{B}_{1, \infty}$
is  defined exactly as  $B_{\kappa}$.
Let $\hat{\cal X}^{c,D}(x)$ denote the total number of the particles
(including the initial one) produced by
the branching process $\hat{B}_{c, D}$ 
starting with a single particle
of  type $x$.

\begin{lem}\label{LO} Let conditions of Theorem \ref{T1A} be fulfilled.
For all large  $D$ and
all small $\varepsilon_1$ in definition of ${\cal A}_n$
one can find $c<1$, arbitrarily  close to
$1$, so that
\[
{\bf P}_{{\cal A}_n}
 \left\{|\tau_n^{U}(x)|<\omega \right\}\leq 
\left(1+b\frac{\log^4n}{n^2}\right)^{n\omega}
 \, {\bf P} \left\{ \hat{\cal X}^{c,D} (x)<\omega
\right\} 
\]
for all large $n$  and all $U \subset \{1, \ldots, n\}$ with $|U|\leq
N\omega$,
where $b$ is some positive constant
independent of $x, c$
and $D$ ($\omega$ and $N$ are defined by (\ref{om}) and (\ref{fN})).
\end{lem}

\noindent
\noindent{\bf Proof.}
At each step of the exploration algorithm which defines $\tau_n^{i, U}$,
the number of the type $y$ offspring 
of a particle of type $x$ has a
binomial distribution $ Bin(N_y',{ p}_{xy}(n))$ where $N_y'$
is the number of remaining vertices of type $y$.

Here we shall explore another 
relation between the binomial and the Poisson  distributions. 
Let $Y_{n,p} \in Bin(n,p)$ and $Z_{\lambda} \in Po(\lambda)$.
Then
it is straightforward to derive from the formulae for the corresponding
probabilities that
for all $0<p<1/4$  and $0\leq k\leq n$ 
\[
{\bf P}
  \{ Y_{n,p}=k\}=\frac{n!}{k!(n-k)!}p^k(1-p)^{n-k}=
  \frac{n!}{n^k(n-k)!}
\left( (1-p)e^{\frac{p}{1-p}}\right)^n e^{-n\frac{p}{1-p}}
  \frac{\left(n\frac{p}{1-p}\right)^k}{k!}
\]
\begin{equation}\label{Re2}
\leq (1+\gamma p^2 )^n \, 
{\bf P}\{ Z_{n\frac{p}{1-p}}=k
  \},
\end{equation}
where $\gamma$ is some positive constant (independent of $n$, $k$ and
$p$). Also notice, that (\ref{Re2}) trivially holds for all $k>n$.

We shall  find now  a lower bound for $N_y'$.
Conditionally on ${\cal A}_n$ we have
\begin{equation}\label{Ma12}
  N_y:=\#\{x_i : x_i=y\} \geq (1-\varepsilon_1)\mu(y)n
\end{equation}
for all $y<D$.
By deleting an  arbitrary set  $U$ with $|U|\leq N\omega$
from $\{1, \ldots, n\}$, we may delete at most $N\omega$ vertices of type $y$.
Hence, conditionally on ${\cal A}_n$  at any step of the 
exploration algorithm which defines $\tau_n^{i,U}$ with
$|\tau_n^{i,U}|<\omega$,
 the number $N_y'$ of the remaining vertices of type $y$, is bounded
from below as follows
\[
N_y' \geq N_y-\omega -N\omega,
\]
and thus
by (\ref{Ma12}) 
\[
N_y' \geq n(1-\varepsilon_1)\mu(y)-\omega -N\omega
\]
for all $y<D$.
Taking into account definitions (\ref{fN}) and (\ref{om}) we derive from here that for any 
$\varepsilon '>0$ one can choose small
$\varepsilon _1>0$  so that
\[
  N_y' \geq (1-\varepsilon')\mu(y)n
\]
for all $y\leq D$ and large $n$. This implies that
conditionally  on ${\cal A}_n$  at any step of the exploration
algorithm we have
\begin{equation}\label{Ma13}
  N_y' 
\frac{{p}_{xy}(n)}{1-{ p}_{xy}(n)} \geq
\mu  (y)(1- \varepsilon')\kappa (x,y)
\end{equation}
for any   $y \leq D$ and large  $n$. 
Now with a help of (\ref{muD}) we rewrite (\ref{Ma13})
as follows: 
\begin{equation}\label{Ma14}
  N_y' 
\frac{{p}_{xy}(n)}{1-{ p}_{xy}(n)} \geq \hat{\mu}_D(y)M_D(1-
\varepsilon')\kappa (x,y) 
=: \hat{\mu}_D(y)c\kappa (x,y) 
\end{equation}
for all  $x,y \in S$, where
\[c=M_D(1- \varepsilon').\]
Recall that $\lim_{D\rightarrow \infty}M_D \uparrow 1$. Therefore
choosing appropriately constants $D$ and $\varepsilon_1$ we can make
$c$ arbitrarily close to $1$.

Using now
relation (\ref{Re2}) between the Poisson and the binomial
distributions, and taking into account  (\ref{Ma14}), we derive for
all $k$ and $N_y'\leq n$
\[
{\bf P}\{ 
Y_{N_y' ,{ p}_{xy}(n)} = k
\}\leq (1+ \gamma{ p}_{xy}^2(n) )^{N_y' } \, 
\, 
{\bf P}
\{ 
Z_{N_y' 
\frac{{p}_{xy}(n)}{1-{ p}_{xy}(n)}}
= k
  \}\]
  \begin{equation}\label{Ma9}
\leq
\left(1+\gamma c_1^2A_1^4\frac{\log^4n}{n^2}
\right)^{n}\, 
{\bf P}
\{ 
Z_{N_y' 
\frac{{p}_{xy}(n)}{1-{ p}_{xy}(n)}}
= k
  \},
\end{equation}
where we used bound (\ref{Kn1}). Note that in the last formula 
$Z_{N_y' 
\frac{{p}_{xy}(n)}{1-{ p}_{xy}(n)}}$ stochastically dominates
$ Z_{\hat{\mu}_D(y)c\kappa (x,y)}$ due to (\ref{Ma14}).
This implies that if conditionally on ${\cal A}_n$,
 at each of at most $\omega$ steps
of the exploration algorithm we replace the
$ Bin(N_y',{ p}_{xy}(n))$ variable 
with the
$$Po\left(\hat{\mu}_D(y)c\kappa (x,y)\right)$$ one, we arrive at
the
following bound using the branching process $\hat{B}_{c, D}$ 
\begin{equation}\label{Ma15}
{\bf P}_{ {\cal A}_n}
 \left\{|\tau_n^{U}(x)|<\omega \right\}
\leq
\left(1+\gamma c_1^2A_1^4\frac{\log^4n}{n^2}
\right)^{n\omega}
 \, {\bf P} \left\{ \hat{\cal X}^{c,D} (x)<\omega
\right\} 
\end{equation}
 for all large $n.$ 
This yields the statement of Lemma \ref{L2}.
\hfill$\Box$

Combining now (\ref{F1}) with (\ref{F4}) and using 
Lemma \ref{LO},
we derive
\[
{\bf P}\left\{  C_1 \Big(
G^{\cal V}(n,\kappa )
    \Big) < \omega \right\}\leq
\left(\left(1+b\frac{\log^4n}{n^2}
\right)^{n\omega}
\sup_{x\in S}{\bf P} \left\{ \hat{\cal X}^{c,D} (x)<\omega 
\right\} 
\right)^{N}
+o(1)\]
\begin{equation}\label{F18}
\leq e^{ b_1\log^3n}
\, \sup_{x\in S}\Big(1-{\bf P} \left\{ \hat{\cal X}^{c,D} (x)\geq\omega 
\right\} 
\Big)^{n/\omega^2} +o(1)
\end{equation}
as $n \rightarrow \infty$, where $b_1$ is some positive constant
independent of $c$
and $D$.

Assume from now on that
\[c=M_D.\]
Define an operator associated with branching process $\hat{B}_{c, D}$:
\[ T_D[f](x):= T_{c\kappa, \hat{\mu}_D}[f](x)=\int_0^D\kappa(x,y)f(y)d\mu(y)
=\sum_{y\leq D}\kappa(x,y)f(y)\mu(y).\]
Clearly, under assumption $\|T_{\kappa}\|<1$ we also have
\begin{equation}\label{A22A}
\|T_{D}\|<1.
 \end{equation}
 Hence, $T_{c\kappa,  \hat{\mu}_D}$ satisfies the conditions of Theorem
\ref{T2} $(ii)$, which together with Remark \ref{R1} implies that
\begin{equation}\label{Re33}
{\hat r}(D):=\sup\{z\geq 1: {\bf E}
z^{\hat{\cal X}^{c,D} (x)} <\infty \} >1
\end{equation}
for all $x\in S$.
It is easy to see that  ${\cal X} (x)$ is stochastically larger than
$\hat{\cal X}^{c,D} (x)$ for all $x\in S$. Therefore we have
\begin{equation}\label{J3}
r_{\kappa}\leq {\hat r}(D)
\end{equation}
for all $D\in S$. Furthermore, we shall prove the following result.

\begin{lem}\label{LO1}
Under the conditions of Theorem \ref{T1A}
\[\lim_{D \rightarrow \infty}{\hat r}( D)=r_{\kappa}.\]
\end{lem}

\noindent
{\bf Proof.} Note that 
 ${\hat r}(D)$ is non-increasing in $D$. Therefore 
 inequality (\ref{J3}) implies existence of the limit
\begin{equation}\label{J4}
\lim_{D \rightarrow \infty}{\hat r}(D)\geq r_{\kappa}.
\end{equation}
We shall show that if 
\begin{equation}\label{J5}
z<{\hat r}( D) \ \ \ \ \ \ \mbox{ for all }  D,
\end{equation}
then also
\begin{equation}\label{J6}
z<{r}_{\kappa}.
\end{equation}
This together with (\ref{J4}) will immediately  imply the statement of the lemma. 

From now on we fix  $z$ which satisfies (\ref{J5}). Then for any $D\in
S $ equation
\begin{equation}\label{J7}
f=ze^{T_D[f-1]}=:\Phi_{D,z}[f]
\end{equation}
has a finite minimal solution $f_D$, which by  (\ref{Jn}) equals
\begin{equation}\label{J8}
f_D(x):=\lim_{k \rightarrow \infty}\Phi_{D,z}^k[1](x)<\infty
\end{equation}
for all $x\in S$.
To prove (\ref{J6}) it is sufficient to show that equation 
\begin{equation}\label{J7*}
f=ze^{T[f-1]}=:\Phi_{z}(f)
\end{equation}
has a finite minimal solution as well.
Therefore we shall prove that
\begin{equation}\label{J8*}
f_{\infty}(x):=\lim_{k \rightarrow \infty}\Phi_{z}^k[1](x)<\infty
\end{equation}
for all $x\in S$, which by Theorem \ref{T3}
is the minimal solution to (\ref{J7*}).

\medskip

\noindent
{\it Claim. } 
For all $k\geq 1$ and for all $x\in S$
\begin{equation}\label{J10}
\lim_{D\rightarrow \infty}
\Phi_{D,z}^k[1](x) = \Phi_{z}^k[1](x).
\end{equation}

\noindent
(Notice that the existence of the limit follows simply by
the monotonicity of  $\Phi_{D,z}$.)
\medskip

\noindent
{\it Proof of the Claim. }
We shall use the induction argument. First, we notice that for all $x\in S$
\[
\Phi_{D,z}[1](x)= z = \Phi_{z}[1](x),
\]
 and 
\begin{equation}\label{J9}
\Phi_{D,z}^2[1](x)= ze^{T_D[1](x)} \uparrow ze^{T[1](x)} = \Phi_{z}^2[1](x)<\infty,
\end{equation}
as $D\rightarrow \infty$.

Assume now that (\ref{J10}) holds
for some $k>1$.
We shall show that then also
\begin{equation}\label{J11}
\lim_{D\rightarrow \infty}\Phi_{D,z}^{k+1}[1](x)= \Phi_{z}^{k+1}[1](x),
\end{equation}
for all $x\in S$, which together with (\ref{J9}) will imply
(\ref{J10})
for all $k\geq 1$.
Set
\[g_D:=\Phi_{D,z}^k[1], \ \ \ \ \ \ \  g:=\Phi_{z}^k[1].\]
By the assumption 
$
g_D \uparrow g $
as $D\rightarrow \infty$.
Then with a help of Theorem on Monotone Convergence we derive
\[\lim_{D\rightarrow \infty}
\Phi_{D,z}^{k+1}[1] =\lim_{D\rightarrow \infty}  \Phi_{D,z}[g_D]=
\lim_{D\rightarrow \infty}
ze^{T_D[g_D-1]} =ze^{T[g-1]}=\Phi_{z}^{k+1}[1], \]
which proves (\ref{J11}). \hfill$\Box$

\medskip

Using (\ref{J10}) we can rewrite function in (\ref{J8*}) as 
\begin{equation}\label{Re30}
f_{\infty}(x)=\lim_{k \rightarrow \infty}
\lim_{D\rightarrow \infty}
\Phi_{D,z}^{k}[1](x).
\end{equation}
Recall that by Theorem \ref{T3} and Remark \ref{R1} we have either
$f_{\infty}(x)<\infty$ or $f_{\infty}(x)=\infty$ for all $x\in S$
(take into account that $S$ is countable here). Our aim is to prove
that $f_{\infty}(x)<\infty$ for all $x\in S$.

Assume, that on the contrary, $f_{\infty}(x)=\infty$ for all $x\in S$.
Let $x_0 = \min \{x\in S\}$ (recall  that $S\subseteq \{1,
2, \ldots\}$).
By (\ref{Re30}) for any $C>0$ there is $k_0=k_0(C)>1$ such that
\begin{equation}\label{Re31}
\lim_{D\rightarrow \infty}
\Phi_{D,z}^{k_0}[1](x_0)>C,
\end{equation}
 which in turn implies that there is $D_0=D_0(C)$ such that
\begin{equation}\label{Re36}
\Phi_{D_0,z}^{k_0}[1](x_0)>C.
\end{equation}
In the case of condition (C3) kernel $\kappa$ is
non-decreasing, and therefore (\ref{Re36}) 
yields as well
\begin{equation}\label{Re31}
\Phi_{D_0,z}^{k_0}[1](x) \geq \Phi_{D_0,z}^{k_0}[1](x_0)>C
\end{equation}
for all $x\in S$.

In the case of condition (C1) (and  (\ref{inf}))
there is a  constant $0<b \leq 1$ such that 
\[\frac{\kappa(x,y)}{ \kappa(x',y)} \geq b\]
for all $x,x', y \in S$. Then for any $D>0$, $k>1$ and for all $x\in S$
\[\Phi_{D,z}^{k}[1](x)=ze^{\sum_S\kappa(x,y)\Big(\Phi_{D,z}^{k-1}[1](y)-1\Big)\mu(y)}
\geq ze^{b
  \sum_S\kappa(x_0,y)\Big(\Phi_{D,z}^{k-1}[1](y)-1\Big)\mu(y)}\geq
\Big(\Phi_{D,z}^{k}[1](x_0)\Big)^b.
\]
This together with (\ref{Re36}) implies 
\begin{equation}\label{Re37}
\Phi_{D_0,z}^{k_0}[1](x) > C^b
\end{equation}
for all $x\in S$.

Now
due to the definition in (\ref{J8}),
and  (\ref{Re37}) or (\ref{Re31}) we have
\begin{equation}\label{Re34}
f_{D_0}(x)=\lim_{k \rightarrow \infty}\Phi_{D_0, z}^k[1](x)
=\lim_{k \rightarrow \infty}\Phi_{D_0, z}^{k}[
\Phi_{D_0,z}^{k_0}[1]
](x)
\geq \lim_{k \rightarrow \infty} \Phi_{D_0,
  z}^{k}[C^b](x).
\end{equation}
It is straightforward to derive taking into account condition (\ref{inf}) and
the definition of $\Phi_{D,
  z}$, that for any $D>0$, $z\geq 1$ and all large $A$ one has
\[ \lim_{k \rightarrow \infty} \Phi_{D,
  z}^{k}[A](x)=\infty.\]
Hence, choosing constant $C$ large enough, we derive from (\ref{Re34})
that 
\[f_{D_0}(x)=\infty, \]
which contradicts inequality in  (\ref{J8}).
Hence, (\ref{J8*}) holds, which finishes the proof of Lemma.

\hfill$\Box$
\medskip

By Lemma \ref{LO1} for any given $\delta_1>0$ we can find a large constant $D$  such that 
\begin{equation}\label{Rev17}
{\hat r}(D)<r_{\kappa}+\delta_1/2.
\end{equation}
From the definition (\ref{Re33}) of ${\hat r}(D)$ it follows that for 
any  $\delta_1>0$ and $x\in S$ there is an unbounded
increasing
sequence $\{\omega_k\}_{k\geq 1}$ of  positive numbers such that
for $c=M_D$
\[
{\bf P}
 \left\{\hat{\cal X}^{c,D}(x)>\omega_k \right\} \geq 
A({\hat r}(D)+\delta_1/2 )^{-\omega _k}\]
for some 
positive constant $A=A(\delta_1,x)<\infty$.
Combining this with (\ref{Rev17}) we get
\begin{equation}\label{Rev16}
{\bf P}
 \left\{\hat{\cal X}^{c,D}(x)>\omega_k \right\} \geq 
A({\hat r}(D)+\delta_1/2 )^{-\omega _k}\geq 
A({r}_{\kappa}+\delta_1 )^{-\omega _k}.
\end{equation}
Note that for any $\delta>0$ and all large $n$ there always exists at
least one $\omega_{k(n)}\in \{\omega_k\}_{k\geq 1}$ such that
\[
\omega =\left(\frac{1}{\log
        r_{\kappa}}-\delta\right) \log n
\leq \omega_{k(n)} \leq
\left(
      \frac{1}{\log
        r_{\kappa}}-\frac{\delta}{2}
  \right)\log n.
\]
  Hence by (\ref{Rev16}) and the fact that $r_{\kappa}>1$ we have
\begin{equation}\label{Rev18}
  {\bf P}
 \left\{\hat{\cal X}^{c,D}(x)>\omega \right\} \geq 
{\bf P}
 \left\{\hat{\cal X}^{c,D}(x)>\omega_{k(n)} \right\} \geq
A({r}_{\kappa}+\delta_1 )^{-\omega_{k(n)} }
\end{equation}
\[\geq
A({r}_{\kappa}+\delta_1 )^{- \left(
      \frac{1}{\log
        r_{\kappa}}-\frac{\delta}{2}
  \right)\log n}.
  \]

Recall that $\kappa$  satisfies  $(C3)$ or $(C1)$
of  Theorem \ref{C1}.
If condition $(C3)$
is satisfied, then 
due to the monotonicity of $\kappa({x,y})$
a vertex of type $x_0=\min S$ has
 among all different types $x\in S$ the smallest
probabilities of the incident edges, which are $\kappa(x_0,y)/n$,
$y\in S$. This together with (\ref{Rev18}) implies for all $x\in S$
\begin{equation}\label{J1}
{\bf P} \left\{ \hat{\cal X}^{c,D} (x)>\omega 
\right\}\geq{\bf P} \left\{ \hat{\cal X}^{c,D} (x_0)>\omega 
\right\}\geq 
A_0({r}_{\kappa}+\delta_1 )^{- \left(
      \frac{1}{\log
        r_{\kappa}}-\frac{\delta}{2}
  \right)\log n}, 
\end{equation}
where $0<A_0=A(\delta_1,x_0)<\infty$.

Otherwise,  when condition $(C1)$ is satisfied, $\kappa(x,y)$ is uniformly bounded from zero and infinity. Then there is a constant $A_0$ such that
$A(\delta_1,x)>A_0>0$ for all $x\in S$, which together with (\ref{Rev18}) yields the same bound 
\[{\bf P} \left\{ \hat{\cal X}^{c,D} (x)>\omega 
\right\}\geq 
A_0({r}_{\kappa}+\delta_1 )^ {- \left(
      \frac{1}{\log
        r_{\kappa}}-\frac{\delta}{2}
  \right)\log n}\]
for all $x\in S$ in this case as well.

Bound (\ref{J1}) allows us to   derive from (\ref{F18}) that for any $\delta >0$ and
$\delta_1>0$
\begin{equation}\label{F17}
{\bf P}\left\{  C_1 \Big(
G^{\cal V}(n,\kappa )
    \Big) < \omega \right\}
\end{equation}
\[\leq e^{ b_1\log^3n}
 \left(1-A_0(r_{\kappa}+\delta_1 )^{-\left(\frac{1}{\log
        r_{\kappa}}-\frac{\delta}{2}\right) \log n }\right)^{
  \frac{n}{ (\alpha\log n)^2} } +o(1),
\]
where $\alpha= \frac{1}{\log
        r_{\kappa}}-\delta$.
Now for  any $\delta >0$ we choose a positive $\delta_1$ so that
\[\gamma_1:= \left(\frac{1}{\log
        r_{\kappa}}-\frac{\delta}{2}\right) \log \left(r_{\kappa}+\delta_1\right) <1.\]
Then (\ref{F17}) becomes
    \begin{equation}\label{F19}
{\bf P}\left\{  C_1 \Big(
G^{\cal V}(n,\kappa )
    \Big)<\left(\frac{1}{\log r_{\kappa}}-\delta\right) \log n\right\}
\leq 
e^{b_1\log^3n} \left(1-\frac{A_0}{ n^{\gamma_1} }\right)^{
  \frac{n}{ (\alpha\log n)^2} } +o(1),
\end{equation}
where the right-hand side goes to zero when $n\rightarrow
\infty$. This
completes the proof of Theorem \ref{Tlb}. \hfill$\Box$

\subsubsection{Proof of Theorem \ref{T1A}.}
Clearly, Theorems \ref{Tlb} and \ref{Tub}  
yield the assertion of Theorem \ref{T1A} when 
$\|T_{\kappa}\|<1$.

When $\|T_{\kappa}\|\geq 1$ we have $r_{\kappa}=1$ by Corollary \ref{C4}.
It is clear that
for any  $0<c < 1/\|T_{\kappa}\|\leq 1$
the size
$ C_1 \Big(
G^{\cal V}(n,\kappa )
    \Big)$ stochastically dominates $ C_1 \Big(
G^{\cal V}(n,c \kappa )
    \Big)$. Then 
we  have
by the previous case for any  $0<c < 1/\|T_{\kappa}\|\leq 1$
\begin{equation}\label{A17}
{\bf P}\left\{ \frac{ C_1 \Big(
G^{\cal V}(n,\kappa )
    \Big)}{\log n} <
   \frac{1}{2 \log r_{c\kappa}} \right\}
\leq {\bf P}\left\{ \frac{ C_1 \Big(
G^{\cal V}(n,c \kappa )
    \Big)}{\log n} <
   \frac{1}{2 \log r_{c\kappa}} \right\}
\rightarrow 0
\end{equation}
as $n\rightarrow \infty$.
By Lemma \ref{LRev} we have $r_{c\kappa}\rightarrow 1$ as $c \uparrow
1/\|T_{\kappa}\|$. Therefore
we derive from (\ref{A17}) that
\[
\frac{C_1 \Big(G^{\cal V}(n,\kappa )
    \Big)}{\log n }
\stackrel{P}{\rightarrow} \infty = \frac{1}{\log r_{\kappa}},
\]
which finishes the proof of Theorem \ref{T1A}. \hfill$\Box$

\bigskip

\bigskip

\noindent
{\bf Acknowledgment} The author thanks P. Kurasov for the
helpful discussions.

\end{document}